\documentclass{amsart}

\let \frak = \mathfrak
\usepackage[mathscr]{euscript}
\let \eusm = \mathscr

\usepackage{amssymb}
\usepackage{amscd}

\newtheorem{thm}{Theorem}[section]

\newtheorem{prop}[thm]{Proposition}
\newtheorem{lem}[thm]{Lemma}
\newtheorem{cor}[thm]{Corollary}
\theoremstyle{definition}
\newtheorem*{defn}{Definition}
\theoremstyle{remark}
\newtheorem*{rem}{Remark}
\newtheorem*{rems}{Remarks}
\newtheorem*{egs}{Examples}

\newcounter{bean}
\newenvironment{mrk}{\begin{list}{\thebean.}{\usecounter{bean}
\setlength{\leftmargin}{0in} \setlength{\itemindent}{.15in}}}{\end{list}}
\newcommand{\R}{\mathbf R}
\newcommand{\C}{\mathbf C}
\newcommand{\E}{\text{End}}
\newcommand{\A}{\text{Aut}}
\renewcommand{\H}{\text{Hom}}
\renewcommand{\r}{\rho}
\newcommand{\ot}{\otimes}
\newcommand{\op}{\oplus}
\renewcommand{\a}{\alpha}
\newcommand{\I}{\text{Ind}}
\renewcommand{\b}{\bar}
\renewcommand{\d}{\delta}
\newcommand{\s}{\sigma}
\renewcommand{\l}{\lambda}
\newcommand{\ov}{\overline}
\newcommand{\D}{\eusm D}
\newcommand{\tD}{\tilde{\eusm D}}
\newcommand{\td}{\tilde}

\title[]{Group actions on central simple algebras}
\author[Daniel S. Sage]{Daniel S. Sage\\
School of Mathematics\\
Institute for Advanced Study\\
Princeton, NJ 08540}
\address{School of Mathematics\\
Institute for Advanced Study\\
Princeton, NJ 08540}

\email{sage@math.ias.edu}
\thanks{Research partially supported by NSF grants DMS 96-233533 and DMS 97-29992.}

\begin{document}

\begin{abstract}
	Let $G$ be a group, $F$ a field, and $A$ a finite-dimensional
central simple algebra over $F$ on which $G$ acts by $F$-algebra
automorphisms.  We study the ideals and subalgebras of $A$ which are
preserved by the group action.  Let $V$ be the unique simple module of $A$.  We show that
$V$ is a projective representation of $G$ and $A\cong\text{End}_D(V)$ makes
$V$ into a projective representation.  We then prove that there is a
natural one-to-one correspondence between $G$-invariant $D$-submodules of
$V$ and invariant left (and right) ideals of $A$.

Under the assumption that $V$ is irreducible, we show that an invariant
(unital) subalgebra must be a simply embedded semisimple subalgebra.  We
introduce induction of $G$-algebras.  We show that each invariant
subalgebras is induced from a simple $H$-algebra for some subgroup $H$ of
finite index and obtain a parametrization of the set of invariant
subalgebras in terms of induction data.  We then describe invariant central
simple subalgebras.  For $F$ algebraically closed, we obtain an entirely
explicit classification of the invariant subalgebras.  Furthermore, we show
that the set of invariant subalgebras is finite if $G$ is a finite group.
Finally, we consider invariant subalgebras when
$V$ is a continuous projective representation of a topological group $G$.
We show that if the connected component of the identity acts
irreducibly on $V$, then all invariant subalgebras are simple.   We then
apply our results to obtain  a particularly nice
solution to the classification problem when $G$ is a compact connected Lie
group and $F=\mathbf C$.
\end{abstract}

\maketitle
\section{Introduction}

Let $G$ be a group, $F$ a field, and $V$ a finite-dimensional F-vector
space on which $G$ acts by $F$-linear automorphisms.  A fundamental problem
in representation theory is to classify the $G$-invariant subspaces of $V$,
in other words, to determine those subspaces of $V$ which inherit a
$G$-action from $V$.  For the case when $G$ is a compact group and
$F=\mathbf C$, this question has been answered completely.  The
representation can be decomposed canonically into a direct sum of
subrepresentations $V=U_1\oplus\dots\oplus U_m$, where each $U_i$ is the
direct sum of $n_i$ copies of an irreducible representation $V_i$ and the
$V_i$'s are pairwise nonisomorphic.  The $G$-invariant subspaces of $U_i$
are parametrized by subspaces of $\mathbf C^{n_i}$ while the
subrepresentations of $V$ are direct sums of subrepresentations of the
$U_i$'s which may be chosen independently. As long as a decomposition of
$V$ into irreducible components is given explicitly (which may be very
difficult in practice), this classification is also entirely explicit.

Let us now replace the vector space $V$ with a finite-dimensional
$F$-algebra $A$.  We suppose further that $A$ is a $G$-algebra, i.e $G$
acts on $A$ by $F$-algebra automorphisms, so that the $G$-action is
well-behaved with respect to ring multiplication.  The natural analogue of
the problem considered above is to determine those $G$-invariant subspaces
of $A$ which have significance in terms of the multiplicative structure of
$A$.  In particular, we would like to classify the $G$-invariant ideals
(left, right, and two-sided) and subalgebras.  These are just special cases
of the general problem of understanding the multiplication of
subrepresentations of $A$.  If $M$ and $N$ are two subrepresentations of
$A$, then $MN$, the $F$-linear span of the set $\{mn\mid m\in M, n\in N\}$,
is also $G$-invariant.  We thus obtain a multiplication on the set of
subrepresentations of $A$.  Invariant ideals and algebras are now easily
expressed in terms of this multiplication; an invariant left ideal is a
subrepresentation $I$ such that $AI\subset I$, an invariant subalgebra is a
subrepresentation $B$ such that $BB\subset B$, and so on.

These problems are much more difficult than the classification of
$G$-invariant subspaces.  It is unreasonable to expect to find a way of
determining $G$-invariant ideals and subalgebras that works for all $A$,
even for $G$ compact and $F=\mathbf C$.  Indeed, if we let $G$ act
trivially on $A$, then this result would give a uniform way of classifying
ideals and subalgebras.  It is thus necessary to limit the class of
algebras under consideration.

In this paper, we restrict attention to central simple algebras over $F$.
Our initial motivation for doing so came from a problem in solid state
physics.  The study of $G$-actions on real and complex central simple
algebras is important in understanding how physical properties such as
conductivity, elasticity, and piezoelectricity of a composite material
depend on the properties of its constituents.  These physical
characteristics are described by elements of a symmetric tensor space
$\text{Sym}^2(T)$, where $T$ is a certain real representation of the
rotation group $SO(n)$.  In general, a property of a composite depends
heavily on the microstructure, i.e. the arrangement of the component
materials.  Let $M\subset\text{Sym}^2(T)$ be the set of all possible values
of a fixed property for composites made with their constituents taken in
prescribed volume fractions.  Typically, $M$ is the closure of an open set
in $\text{Sym}^2(T)$ and may be described by a system of inequalities, so
that away from the boundary of $M$, it is possible to make any desired
small change in the property by varying the microgeometry.  However, in
certain unusual situations, some of the inequalities become equations,
determining a proper closed submanifold $E$ in which $M$ is locally closed.
The submanifold $E$ and also the equations defining $E$ are called exact
relations for the property.  The variability of the property with
microstructure is thus drastically reduced when an exact relation is
present.  Recent work of Grabovsky, Milton, and Sage has shown how to
classify exact relations in terms of the multiplication of
$SO(n)$-subrepresentations in the endomorphism algebra $\text{End}_\R(T)$;
in particular, invariant algebras and ideals of this central simple algebra
have great physical significance \cite{G,GS,GMS}.

Let $A$ be a central simple algebra over $F$, and suppose $G$ acts on $A$
by $F$-algebra automorphisms.  In the first part of this paper, we show
that $A$ is isomorphic to the algebra of $D$-endomorphisms of a
projective representation $V$ of $G$, where $D$ is a certain central
division algebra.  We then prove that there is a natural one-to-one
correspondence between $G$-invariant $D$-submodules of $V$ and invariant
left (and right) ideals of $A$.  Indeed, we show that if $G$ is compact and
$A$ is the endomorphism algebra of a complex representation, then the
parametrization of invariant left and right ideals of $A$ is the same as
the classical parametrization of invariant subspaces of $V$ described
above.  In particular, $V$ is irreducible if and only if there are no
invariant proper left (right) invariant ideals, and $V$ is multiplicity
free  if and only if there are a finite number of left
(right) invariant ideals.

In the second part of the paper, we turn to the much more complicated
problem of understanding unital invariant subalgebras of $A$ under the
additional hypothesis that $V$ is an irreducible projective representation.
We show that an invariant subalgebra $B$ must be a simply embedded
semisimple subalgebra; this means that both $B$ and its centralizer in $A$
must be direct products of isomorphic simple algebras.  We then introduce
induction of $G$-algebras.  We show that each invariant subalgebras is
induced from a simple $H$-algebra for some subgroup $H$ of finite index and
obtain a parametrization of the set of invariant subalgebras in terms of
induction data.  We then describe invariant central simple subalgebras.

Combining these two results, we obtain an entirely explicit classification
of the invariant subalgebras for $F$ algebraically closed.  This
classification shows that the set of invariant subalgebras of $A$ encodes
complicated information about $G$ and $V$, involving both how $V$ can be
expressed as an induced representation $\I_H^G(W)$ and how $W$ can be
factored into the tensor product of projective representations.  It should
be observed that for $F=\C$ and $G$ finite, knowing the character table of
$G$ does not suffice to determine all invariant subalgebras.  In fact, even
in the simplest case where $V$ is a primitive representation, the character
table of a covering group of $G$ is needed to find all invariant
subalgebras.  When $V$ is primitive, we show that the only nonunital
invariant subalgebra is $\{0\}$.  Finally, we prove that for $G$ finite and
$F$ algebraically closed, the set of invariant subalgebras is finite,
and we describe how finiteness fails in the general case.

In the final section of the paper, we consider invariant subalgebras when
$V$ is a continuous projective representation of a topological group $G$.
We show that if the connected component of the identity acts
irreducibly on $V$, then all invariant subalgebras are simple.  We then
apply our results to obtain a theorem of Etingof giving a particularly elegant
solution to the classification problem when $G$ is a compact connected Lie
group and $F=\C$.  In fact, suppose $G$ is semisimple and simply connected,
say $G=G_1\times\dots\times G_n$ with each $G_i$ simple.  The
representation $V$ is then isomorphic to $V_1\otimes\dots\otimes V_n$, for
some irreducible representations $V_i$ of $G_i$.  We show that the
$G$-invariant subalgebras of $A$ are parametrized by the subsets $J$ of
$\{i\mid V_i\ne\C\}$ via $J\mapsto \bigotimes_{j\in J}\text{End}_{\mathbf
C}(V_j)$ and that the only nonunital invariant subalgebra is $\{0\}$.  In
particular, if $G$ is simple, the invariant subalgebras are $\mathbf C$ and
$A$.

We have also obtained results on the general problem of multiplication of
subrepresentations in central simple algebras when $G$ is a compact, simply
reducible group.  This means that $g$ and $g^{-1}$ are conjugate for all
$g\in G$ (so that all $G$-modules are self-dual) and $G$ is
multiplicity-free, i.e. if $V$ and $W$ are irreducible, then each isotypic
component of $V\otimes W$ is irreducible.  (The most familiar examples of
simply reducible groups are $S_3$, $S_4$, the quaternion group, $SU(2)$,
and $SO(3)$.)  However, since the proofs use quite different techniques,
these results will appear in another paper \cite{S}.

It is a great pleasure to thank Yury Grabovsky for first bringing these
problems to my attention and for explaining their importance in physics.  I
would also like to thank Daniel Allcock for several helpful comments and
Pavel Etingof for letting me use his unpublished result on invariant
subalgebras of compact connected Lie groups.

\section{Preliminaries and Invariant Ideals}

Let $A$ be a finite-dimensional central simple algebra over the field $F$,
and let $V$ be a simple (left) $A$-module.  The module $V$ is unique up to
isomorphism and is a finite-dimensional vector
space over $F$.  By Schur's Lemma, the ring $D=\E_A(V)$ is a central
division algebra, and $V$ is naturally a left $D$-module.  It is well-known
that $A$ is isomorphic to $\E_D(V)$, and from now on, we assume
without loss of generality that $A=\E_D(V)$.

It is easy to construct examples of central simple algebras on which the
group $G$ acts by $F$-algebra automorphisms.  Recall that a mapping
$\rho:G\to GL(V)$ is called a projective representation of $G$ over $F$ if
$\rho(1)=1_V$ and if there exists $\alpha:G\times G\to F^*$ such that
$\rho(xy)=\alpha(x,y)\rho(x)\rho(y)$ for all $x,y\in G$.  (Equivalently, we
can view a projective representation as a homomorphism $G\to PGL(V)$.)  The
map $\a$ is a $2$-cocycle.  Let $\b{g}$ be the basis vector corresponding
to $g\in G$ in the twisted group algebra $F^{\a}G$.  A projective
$\a$-representation is just an $F^{\a}G$-module via $\b{g}v=\rho(g)(v)$, and
we also use this notation.  (For linear representations, we just write
$gv$.)  The map $\pi:G\to GL(A)$ then makes $A$ into a (linear)
representation of $G$ with $(\pi(g)f)(v)=\rho(g)(f(\rho(g)^{-1}(v))$ for
all $g\in G$, $f\in A$, and $v\in V$.  Moreover, the linear map $\pi(g)$ is
in fact an algebra automorphism.  It turns out that all central simple
algebras on which $G$ acts via algebra automorphisms are of this type.

\begin{prop}  \label{T:SN} Suppose that $G$ acts on $A=\E_D(V)$ by $F$-algebra
automorphisms, i.e. $A$ is a representation of $G$ via a homomorphism
$G\overset{\pi}{\to}\A(A)$.  Then $V$ is a projective representation of
$G$ determined up to projective equivalence, and the $G$-action on $A$ is the natural action induced by the
projective $G$-action on $V$.
\end{prop}
\begin{proof}  Any automorphism of $A$ is inner by the Skolem-Noether
theorem.  Hence, we obtain a function $\hat{\rho}:G\to A^{\times}\subset
GL(V)$ such that $\pi(g)(a)=\hat{\rho}(g)a\hat{\rho}(g)^{-1}$ for all
$g\in G$ and $a\in A$.  Since $\pi(1)=1_A$, we have $\hat{\rho}(1)\in
Z(A)^{\times}=F^*$.  Setting $\rho(g)=\hat{\rho}(g)/\hat{\rho}(1)$ gives
$\rho(1)=1_V$.  Also, the equation $\pi(gh)=\pi(g)\pi(h)$ implies that
$\rho(gh)\rho(h)^{-1}\rho(g)^{-1}$ is central and therefore a nonzero
multiple of the identity.  It follows that $(V,\rho)$ is a projective
representation of $G$ giving rise to $\pi$.
\end{proof}
 
We now recall the ideal structure of $A$.  Let $\eusm S(V)$ denote the set
of $D$-subspaces of $V$ partially ordered by inclusion.  This poset is in
fact a complete lattice, with the greatest lower bound and least upper
bound of a collection of subspaces given by their intersection and sum
respectively.  Similarly, the sets $\eusm L(A)$ and $\eusm R(A)$ of left
and right ideals of $A$ are complete lattices.  It will be convenient to
work with the dual lattice $\eusm L(A)^*$ of left ideals under reverse
inclusion (and with the supremum and infimum reversed).  If $L$ is a
$D$-submodule of $V$, we define the annihilator and coannihilator of $L$
by $\text{Ann}(L)=\{f\in A\mid f(L)=0\}$ and $\text{Coann}(L)=\{f\in A\mid
f(V)\subset L\}$; these are respectively left and right ideals of $A$.  We
then have the well-known fact that all left and right ideals of $A$ are of
this form.
\begin{prop} The maps $\eusm S(V)\overset{\text{Ann}}{\to}\eusm L(A)^*$
and $\eusm S(V)\overset{\text{Coann}}{\longrightarrow}\eusm R(A)$ are isomorphisms of
complete lattices.  The inverses are given by $I\mapsto\bigcap_{f\in
I}\text{Ker}(f)$ and $J\mapsto\sum_{f\in J}f(V)$, where $I\in\eusm L(A)$
and $J\in\eusm R(A)$.
\end{prop}
\begin{rem}  In matrix language, this simply says that a left ideal
consists of all matrices (with respect to some basis depending on the
ideal) with zeroes in given columns while a right ideal consists of
all matrices with zeros in given rows.
\end{rem}

Let $\eusm S_G(V)\subset\eusm S(V)$ be the complete sublattice of all
$D$-subspaces of $V$ preserved by the $G$-action on $V$.  Similarly,
we define the complete sublattices $\eusm L_G(A)\subset\eusm L(A)$ and
$\eusm R_G(A)\subset\eusm R(A)$ of $G$-invariant left and right ideals of
$A$.  It is natural to conjecture that the sublattices $\eusm L_G(A)$ and
$\eusm R_G(A)$ are just the images of $\eusm S_G(V)$ under the above
isomorphisms, i.e. invariant left and right ideals are annihilators and
coannihilators respectively of subrepresentations of $V$.  This is indeed
the case.
\begin{thm}\label{T:ideal}  The restrictions of the maps $\text{Ann}$ and
$\text{Coann}$ define isomorphisms of complete lattices $\eusm
S_G(V)\overset{\text{Ann}}{\to}\eusm L_G(A)^*$ and $\eusm
S_G(V)\overset{\text{Coann}}{\longrightarrow}\eusm R_G(A)$.
\end{thm}
\begin{proof}  In order to prove the first isomorphism, it suffices to
show that $\text{Ann}(\eusm S_G(V))\subset\eusm L_G(A)^*$ and
$\text{Ann}^{-1}(\eusm L_G(A)^*)\subset\eusm S_G(V)$.  If $L$ is a
subrepresentation of $V$ and $f\in\text{Ann}(L)$, then $(g\cdot
f)(v)=\b{g}(f(\b{g}^{-1}(v)))=\b{g}(0)=0$ for all $g\in G$ and $v\in
L$.  Thus, $\text{Ann}(L)$ is $G$-invariant.  Conversely, if $I$ is an
invariant left ideal and $v\in\text{Ann}^{-1}(I)=\bigcap_{f\in
I}\text{Ker}(f)$, then we also have $v\in\bigcap_{f\in
I}\text{Ker}(g\cdot f)$.  Since $\rho(g)$ is bijective, this gives
$f(\b{g}^{-1}v)=0$ for all $g\in G$ and $f\in I$.  It follows that
$\text{Ann}^{-1}(I)$ is $G$-invariant.

The proof for invariant right ideals is similar.
\end{proof}
\begin{rems}  \begin{mrk} \item Since $A$ is simple, the only two-sided ideals are $\{0\}$ and $A$
which are of course $G$-invariant.  However, it is a general fact that if
$B$ is an arbitrary $G$-algebra on which $G$ acts by inner automorphisms,
then all two-sided ideals are $G$-invariant.  Indeed, if $I$ is a two-sided
ideal and the action of $g$ on $B$ is given by conjugation by $b_g\in
B^{\times}$, then $gI=b_g I b_g^{-1}\subset I$.
\item Suppose that $F$ is algebraically closed and $V$ is a completely
reducible linear representation of $G$, say  $V\cong n_1V_1\oplus\dots\oplus n_mV_m$ where the $V_i$'s are pairwise
nonisomorphic irreducible representations.  Then the $G$-invariant left
(and right) ideals of $\E_F(V)$ are parametrized by
$\prod_{i=1}^m\{\text{subspaces of }F^{n_i}\}$. \end{mrk}
\end{rems}

This theorem allows us to characterize certain properties of
representations in terms of the associated endomorphism algebras.
\begin{cor}\begin{enumerate} \item The projective representation $V$ is
irreducible if and only if $\E_F(V)$ has no proper invariant one-sided
ideals. 
\item Let $D$ be a central division algebra, and suppose $V$ is a
$D$-module on which $G$ acts (projectively) by $D$-linear automorphisms.
Then $V$ is $D$-irreducible (i.e. has no $G$-invariant $D$-submodules) if
and only if $\E_D(V)$ has no proper invariant one-sided
ideals.
\item Suppose that $F$ is an infinite field and $V$ is  completely reducible.  Then $V$ is multiplicity free if and only if
$\E_F(V)$ has a finite number of  invariant one-sided
ideals.
	   \end{enumerate}
\end{cor}
\begin{proof} The first two statements are clear from the theorem.  The
last follows from the second remark and the fact that for an infinite
field, a vector space has an infinite number of subspaces if and only if it
has dimension larger than one.
\end{proof}
It is worth noting that in spite of the strong connection between
subrepresentations and invariant ideals, the group action on a
subrepresentation does not determine the action on the corresponding left
and right invariant ideals or vice versa.

If $B$ is a semisimple (finite-dimensional) algebra on which $G$ acts by inner
automorphisms, this theorem can be used to determine the invariant ideals
of $B$.  Let $B=B_1\oplus\dots\oplus B_s$ where the simple component $B_i$ can be viewed as
$\E_{D_i}(V_i)$ where $D_i$ is a finite-dimensional division algebra over
$F$ and $V_i$ is a finite-dimensional $D_i$-module.  By the first remark, the
two-sided ideal $B_i$ is invariant and is thus a simple algebra on which
$G$ acts by inner automorphisms.  Since left and right ideals of $B$ are
just direct sums of left and right ideals of $B_i$, we obtain the following
corollary:

\begin{cor} The maps $\prod_{i=1}^s \eusm S_G(V_i)\to\eusm L_G(B)^*$ and
$\prod_{i=1}^s \eusm S_G(V_i)\to\eusm R_G(B)$ given
by $(L_1,\dots,L_s)\mapsto \bigoplus_{i=1}^s\text{Ann}(L_i)$ and
$(L_1,\dots,L_s)\mapsto \bigoplus_{i=1}^s\text{Coann}(L_i)$ respectively
are isomorphisms of complete lattices.
\end{cor}

It is also possible to obtain analogous results for certain spaces
of homomorphisms between two representations of $G$.  Let $V$ and $W$ be
two finite-dimensional linear representations over the division algebra $D$ (which as
above is finite-dimensional over $F$) on which $G$ acts by $D$-linear
automorphisms.  The $F$-vector space $\H_D(V,W)$ is a representation of $G$
via the action $(g\cdot f)(v)=g(f(g^{-1}v))$; moreover, it has the
structure of an $(\E_D(W),\E_D(V))$-bimodule.  Let $L$ be a $D$-submodule
of $V$, and define the annihilator of $L$ by $\text{Ann}(L)=\{f\in
\H_D(V,W)\mid f(L)=0\}$.  This is a (left) $\E_D(W)$-submodule, and all
$\E_D(W)$-submodules of $\H_D(V,W)$ are of this form.  In fact, if we
denote the lattice of $\E_D(W)$-submodules of $\H_D(V,W)$ by $\eusm
L(\H_D(V,W))$, the map $\eusm S(V)\overset{\text{Ann}}{\to}\eusm
L(\H_D(V,W))^*$ is an isomorphism of complete lattices.  Similarly, for $M$
a $D$-submodule of $W$, we have the (right) $\E_D(V)$-submodule
$\text{Coann}(M)=\{f\in \H_D(V,W)\mid f(V)\subset L\}$, and denoting the
lattice of $\E_D(V)$-submodules of $\H_D(V,W)$ by $\eusm R(\H_D(V,W))$, we
obtain the isomorphism of complete lattices $\eusm
S(W)\overset{\text{Coann}}{\longrightarrow}\eusm R(\H_D(V,W))$.  Again, we can characterize the sublattices $\eusm L_G(\H_D(V,W))$ and $\eusm
R_G(\H_D(V,W))$ of $G$-invariant submodules in terms of these isomorphisms.

\begin{thm}  The restrictions of the maps $\text{Ann}$ and
$\text{Coann}$ define isomorphisms of complete lattices $\eusm
S_G(V)\overset{\text{Ann}}{\to}\eusm L_G(\H_D(V,W))^*$ and $\eusm
S_G(W)\overset{\text{Coann}}{\longrightarrow}\eusm R_G(\H_D(V,W))$.
\end{thm}
\begin{proof} The proof is similar to the proof of Theorem \ref{T:ideal}.
\end{proof}

\begin{rem} The only sub-bimodules of $\H_D(V,W)$ are $\{0\}$ and $\H_D(V,W)$,
which are $G$-invariant.
\end{rem}

\section{Invariant Subalgebras}

We now turn our attention to subalgebras of the algebra $A=\E_D(V)$ which
are preserved by the group action.  (All subalgebras will be assumed
to contain $1$ unless otherwise specified.)  Invariant subalgebras are much more
difficult to understand than invariant ideals, and in general, invariant
subalgebras can be very badly behaved. For example, if we let $G$ act
trivially on $\E_F(V)$, then every subalgebra is invariant.  This means
that if $V$ has dimension $n$, then $\E_F(V)$ contains every
$n$-dimensional $F$-algebra as an invariant subalgebra.  Moreover, it is
not even true that the ring of invariants $A^G$ need be semisimple, if $G$
is infinite or $G$ is finite with the characteristic of $F$ dividing $|G|$
\cite{M}.  We will therefore need to place additional restrictions on the
$G$-algebra $A$.

\begin{defn}  A central simple $G$-algebra over $F$ is called {\em
$G$-simple} if the associated projective representation $V$ is
irreducible.
\end{defn}

We assume from now on that $A$ is $G$-simple.  Note that under this
hypothesis, the possible pathologies involving $A^G$ are avoided, since by
Schur's lemma, $A^G$ is a division algebra.

We will show that all $G$-invariant subalgebras of $A$ are semisimple with
a very special structure.  Indeed, we will give a complete classification
when $F$ is algebraically closed.  As a first step, we have the following
result:
\begin{prop} \label{T:ss} Let $B$ be an invariant subalgebra of $A$.  Then $B$ is
semisimple, and the Wedderburn components of $B$ are all isomorphic as
$F$-algebras.  Moreover, if $U$ is any simple $B$-submodule of $V$, then
for each $g\in G$, $\b{g}U$ is also a simple $B$-submodule, and  any
simple $B$-module is isomorphic to some $\b{g}U$.
\end{prop}
\begin{proof}  The inclusion of $B$ in $A$ makes the $A$-module $V$ into a
$B$-module.  Let $U$ be a simple $B$-submodule of $V$; for example, take
$U$ to be a $B$-submodule of minimal dimension as an $F$-vector space.
Consider the translate $\b{g}U$ for $g\in G$. Note that the
$G$-invariance of $B$ implies that
\begin{equation}\label{E:gU}  b\b{g}(u)=\b{g}\b{g}^{-1}b\b{g}(u)=\b{g}(g^{-1}\cdot
b)(u)\in \b{g}U
\end{equation}
 for all $b\in B$ and $u\in U$.  Here, we have used the fact that
$\ov{g^{-1}}=\a(g,g^{-1})\b{g}^{-1}$, where $\a$ is the cocycle defined by
$(V,\rho)$. Thus, $\b{g}U$ is a $B$-submodule of $V$.  Moreover, $\b{g}U$ is
simple, since the same argument shows that if $W$ is a submodule of
$\b{g}U$, then $\ov{g^{-1}}W$ is a submodule of $U$.  The sum $\sum_{g\in
G}\b{g}U$ is evidently a nonzero $G$-invariant subspace of $V$, and by
irreducibility, $V=\sum_{g\in G}\b{g}U$.  Thus, $V$ is a semisimple
$B$-module, and we can choose $g_1,\dots,g_l\in G$ such that
$V=\oplus_{i=1}^l \b{g_i}U$.

Let $u_1,\dots, u_k$ be an $F$-basis for $U$.  The map $B\to
\oplus_{i=1}^l k(\ov{g_i}U)$ given by
$b\mapsto(b\ov{g_1}u_1,\dots,b\ov{g_1}u_k,\dots,b\ov{g_l}u_1,\dots,b\ov{g_l}u_k)$ is a
$B$-homomorphism.  If $b$ is in the kernel, then $b$ kills an $F$-basis of
$V$, and since $b\in A\subseteq\E_F(V)$, we have $b=0$; hence, the map is
injective.  This shows that $B$ is a semisimple $F$-algebra, and any simple
$B$-module is isomorphic to $\b{g}U$ for some $g\in G$.  The simple components
of $B$ are of the form $\E_{D_g}(\b{g}U)$, where $D_g=\E_B(\b{g}U)$.  To complete
the proof, it suffices to verify that $\E_{D_g}(\b{g}U)$ is isomorphic to
$\E_{D'}(U)$, where $D'=D_1$.

We first show that the division algebras $D'$ and $D_g$ are isomorphic via
the map $d\to\b{g}d\b{g}^{-1}$.  Using the formula for the $B$-action on
$\b{g}U$ given in \eqref{E:gU}, we have
$\b{g}d\b{g}^{-1}(b\b{g}u)=\b{g}d\b{g}^{-1}(\b{g}(g^{-1}\cdot
b)u))=\b{g}d((g^{-1}\cdot b)u)=\b{g}(g^{-1}\cdot
b)d(u)=b\b{g}d(u)=b\b{g}d\b{g}^{-1}(\b{g}u)$ for all $d\in D'$ and $u\in U$, so that $\b{g}d\b{g}^{-1}\in
D_g$.  It is clear that this is an $F$-algebra homomorphism.  In fact, it is
an isomorphism with inverse map $D_g\to D'$ given by $\hat
d\mapsto\ov{g^{-1}}\hat d\ov{g^{-1}}^{-1}$.  This follows since $\b{g}\ov{g^{-1}}\in
F^*$ and elements of $D'$ and $D_g$ are $F$-linear.

Now suppose $f\in\E_{D'}(U)$.  The $F$-map $\b{g}f\b{g}^{-1}:\b{g}U\to \b{g}U$ is
$D_g$ linear as$\b{g}f\b{g}^{-1}(\b{g}d\b{g}^{-1}(\b{g}u))=\b{g}f(du)=\b{g}df(u)=
\b{g}d\b{g}^{-1}(\b{g}f\b{g}^{-1}(\b{g}u))$. Thus, we have an $F$-algebra
homomorphism $\E_{D'}(U)\to\E_{D_g}(\b{g}U)$, $f\mapsto\b{g}f\b{g}^{-1}$, which is
in fact an isomorphism with inverse $\hat f\mapsto\ov{g^{-1}}\hat
f\ov{g^{-1}}^{-1}$.
\end{proof}

\begin{cor} The invariant subalgebra $B$ is simple if and only if any
for any simple $B$-submodule $U$ of $V$, the $B$-modules $U$ and $\b{g}U$ are
isomorphic for all $g\in G$.
\end{cor}

Although the proposition places significant restrictions on the structure
of a $G$-invariant subalgebra, it turns out that the subalgebra must
satisfy a much more stringent condition which depends on the ambient
algebra $A$.  For the time being, let $B=B_1\op\dots\op B_l$ be an
arbitrary semisimple subalgebra of $A=\E_D(V)$ where the $B_i$'s are simple
$F$-algebras with corresponding simple modules $W'_i$.  Note that
$V=\oplus_{i=1}^l m'_iW'_i$ with positive multiplicities $m'_i$ (or else
$1_{B_i}(V)=0$ for some $i$, contradicting the fact that the central
primitive idempotent $1_{B_i}$ is a nonzero element of $A$).  The
subalgebra $B$ consists of $D$-linear maps, so $V$ can also be viewed as a
$(B,D^{op})$-bimodule or equivalently as a $B\otimes_F D$-module.  Since
$D$ is central simple, $B\otimes D$ is semisimple with Wedderburn
components $B_1\ot D,\dots, B_l\ot D$, and we can write
$V=\oplus_{i=1}^l m_iW_i$, where the $W_i$'s are the simple $B\otimes
D$-modules, with $W_i=W'_i\ot D$ isomorphic to a minimal left ideal of $B_i\ot_F D$.
Again, each $m_i$ is nonzero.  In fact, we can say more.

\begin{lem}  Let $V'_i$ and $V_i$ be the isotypic $B$ and $B\ot_F
D$-submodules of $V$ for $W'_i$ and $W_i$ respectively.  Then $V'_i=V_i$
and $m_i=m'_i/\dim_F D$.  In particular, $V'_i$ is a $D$-submodule of $V$.
\end{lem}
\begin{proof}  Recall that $V'_i$ and $V_i$ are
the one-eigenspaces of the central primitive idempotents $1_{B_i}$ and
$1_{B_i}\ot 1_D$.  Since these are the same maps on $V$, we have $V'_i=V_i$
for all $i$.  Also $\dim W_i=(\dim D)(\dim W'_i)$, so $m'_i=m_i\dim D$.
\end{proof}

\begin{defn}  A semisimple subalgebra $B$ of $A=\E_D(V)$ is called {\em
symmetrically embedded} if the Wedderburn components of $B$  are all isomorphic as
$F$-algebras and if the simple $B$ modules appearing in $V$
have the same multiplicity $m'$, i.e. if $m'=m'_1=\dots =m'_l$.
(It is equivalent to replace either condition with the analogous statement
involving $B\ot D$.)
\end{defn}

More explicitly, an element $b\in B$ acts on each copy of $W_i$ in the same
way, so is represented by a block diagonal matrix (in $M_{\dim_D
V}(D^{op})$) with $m_1+\dots+m_l$ blocks, $m_i$ of which consist of the
$\dim_D W_i\times\dim_D W_i$ matrix corresponding to $b|_{W_i}$.  If $B$ is
symmetrically embedded, then the blocks are all the same size and for each
$i$, the matrix for $b|_{W_i}$ appears $m$ times.  Note that this implies
that $\dim_D V=ml\dim_D W_i=ml\dim_F W'_i$ for any $i$.

It is clear from the above lemma that whether a subalgebra satisfies the
above property does not depend on the central division algebra $D$.
Indeed, we have:
\begin{prop}  Suppose that $V$ is a module for two central division
algebras $D$ and $D'$.  If $B$ is a semisimple subalgebra of both $A=\E_D(V)$
and $A'=\E_{D'}(V)$, then $B$ is symmetrically embedded in $A$ if and only
if it is symmetically embedded in $A'$.
\end{prop}
\begin{proof} Since both $A$ and $A'$ are subalgebras of $\E_F(V)$, we can
assume without loss of generality that $D'=F$, and this case follows
immediately from the lemma. 
\end{proof}

In order to see the importance of symmetrically embedded subalgebras, we
need to recall some information about centralizers of semisimple
subalgebras of central simple algebras.  Let $Z_A(B)$ denote the
centralizer in $A$ of the subalgebra $B$.  We call $B$ a {\em Howe
subalgebra} if it equals its double centralizer $Z_AZ_A(B)$ and say that
the pair $(B,Z_A(B))$ is a {\em dual pair}.  A strong version of the Double
Centralizer Theorem  states that if $B$ is
semisimple, then $Z_A(B)$ is also semisimple and $B$ is a Howe
subalgebra \cite[Theorem 4.10]{J}. In other words, the mapping $B\mapsto Z_A(B)$ provides a
duality operator on the set of semisimple subalgebras of $A$.  It is
possible to calculate the Wedderburn structure of $Z_A(B)$ by an argument
due to Moeglin, Vign{\'e}ras, and Waldspurger \cite[p.12]{MVW}.  Note that
$f\in \E_{B\ot_F D}(V)$ if and only if $f$ is a $D$-linear map which
commutes with the action of $B$, i.e. if and only if $f\in Z_A(B)$.  Using
the decomposition $V=\op_{i=1}^l m_iW_i$ of $V$ into simple $B\ot_F
D$-submodules, it is immediate that $\E_{B\ot D}(V)\cong
\oplus_{i=1}^l\E_{B\ot D}(m_iW_i)\cong\op_{i=1}^l M_{m_i}(D_i)$, where
$D_i=\E_{B\ot_F D}(W_i)$ is a division algebra over $F$.  Since
$W_i=W'_i\ot D$, $D_i$ is canonically isomorphic to $\E_B(W'_i)$.  Summing up, we
have:
\begin{thm} \label{T:DCT} There is a duality on the set of semisimple  subalgebras of
$A$ given by $B\mapsto Z_A(B)$ which
preserves the number of Wedderburn components of the subalgebras.  Moreover,
if $V\cong\op_{i=1}^l m_iW_i$ is the decomposition of $V$ into simple
$B\ot_F D$-modules and $D_i$ is the division algebra $\E_B(W'_i)=\E_{B\ot_F D}(W_i)$,
then $Z_A(B)\cong \op_{i=1}^l M_{m_i}(D_i)$.
\end{thm}

There are also two maps from a semisimple subalgebra to the set of
self-dual (i.e. commutative) semisimple subalgebras, given by $B\mapsto
Z(B)$, the center of $B$, and $B\mapsto Z_0(B)$, the $F$-linear span of the
central primitive idempotents of $B$.  These are respectively the largest
and smallest self-dual subalgebras with the same central primitive
idempotents as $B$.  It
is clear that both maps are constant on dual pairs.  With the notation of
the theorem, $Z(B)\cong\op_{i=1}^l Z(D_i)$ and $Z_0(B)\cong F^l$.

We can now reformulate the concept of a symmetrically embedded subalgebra
in terms of centralizers.
\begin{prop}  \label{T:embed} A semisimple subalgebra $B$ is symmetrically embedded in $A$
if and only if both $B$ and $Z_A(B)$ are direct sums of isomorphic simple
$F$-algebras.  
\end{prop}
\begin{proof}  Suppose $B$ is symmetrically embedded.  By definition, the
Wedderburn components of $B$ are all isomorphic.  The Jacobson density
theorem implies that $B\ot_F D\cong \E_{D_i}(W_i)$ for all $i$, and by the
structure theorem for simple Artinian algebras, the $D_i$'s are all
isomorphic as $F$-algebras.  Since the multiplicities of the $W_i$'s in $V$
are the same, it follows from Theorem \ref{T:DCT} that $Z_A(B)$ is a direct
sum of isomorphic simple $F$-algebras.

Conversely, suppose that both $B$ and $Z_A(B)$ have isomorphic Wedderburn
components.  Then $M_{m_i}(D_i)\cong M_{m_j}(D_j)$ for all $i$ and $j$, and
so the $m_i$'s are equal by the structure theorem for simple Artinian
algebras.  Thus, $B$ is symmetrically embedded in $A$.
\end{proof}

Next, we need an easy, but important lemma on centralizers of invariant
subalgebras.
\begin{lem}\label{T:cent}  Let $R$ be a $G$-algebra, and $S$ a $G$-invariant subalgebra.
Then the centralizer $Z_R(S)$ is also an invariant subalgebra.  In
particular, the center of $S$ $Z(S)=Z_S(S)$ is an invariant subalgebra.
\end{lem}
\begin{proof}  This follows immediately from the fact that $(g\cdot
z)s=g\cdot(z(g^{-1}\cdot s))=g\cdot((g^{-1}\cdot s)z)=s(g\cdot z)$ for all
  $g\in G$, $s\in S$, and $z\in Z_R(S)$.
\end{proof}

Combining the lemma with Propositions \ref{T:ss} and \ref{T:embed}, we
obtain the theorem:
\begin{thm} Let $B$ be an invariant subalgebra of $A$.   Then $B$ is
symmetrically embedded in $A$.
\end{thm}

We now describe a fundamental construction of invariant subalgebras.  We
will then show that all invariant subalgebras are of this type and obtain a
classification of them.

We first need to introduce induction of $G$-algebras.  Let $H$ be a
subgroup of finite index in $G$, and suppose that $C$ is an $H$-algebra.
We show how to define a natural $G$-algebra structure on $\I_H^G(C)$ making
$\I_H^G$ into a functor from the category of $H$-algebras into the category
of $G$-algebras.
\begin{prop} \label{T:algind} There is a unique $G$-algebra structure on
$\I_H^G(C)=FG\otimes_{FH}C$ extending the $H$-algebra $1\otimes C$ such
that distinct $G$-translates of $1\ot C$ annihilate each other.  If
$\{g_1,\dots,g_n\}$ is a left transversal for $H$ in $G$, then the algebra
multiplication is given by $(g_i\ot b)(g_j\ot b')=\delta_{ij}(g_i\ot bb')$
for $b,b'\in C$.  As $F$-algebras, $\I_H^G(C)$ is isomorphic to $C^n$.
Furthermore, this definition makes $\I_H^G$ into a functor from the category
of $H$-algebras into the category of $G$-algebras.
\end{prop}
\begin{proof}  Uniqueness is clear.  To show existence,
recall that the coinduced representation $\H_{FH}(FG,C)$ (with $G$ acting
by $(g\cdot f)(x)=f(xg)$ for $x,g\in G$) is isomorphic to
$\I_H^G(C)$ via the map $\phi\mapsto\sum_{i=1}^n g_i\ot\phi(g_i^{-1})$.  If
$\phi$ and $\psi$ are $FH$-linear, then $\phi\psi$ is as well, since
$(\phi\psi)(hy)=(h\phi(y))(h\psi(y))=h\cdot(\phi\psi(y))$ for $h\in H$ and
$y\in FG$.  Thus, pointwise multiplication makes $\H_{FH}(FG,C)$ into a
$G$-algebra; translating the multiplication back to $\I^G_H(C)$ gives the
desired formula.    The elements $g_i\ot 1$ are pairwise orthogonal central
idempotents summing to the identity element in $\I_H^G(C)$, which is thereby
isomorphic to $\bigoplus_{i=1}^n (g_i\ot C)\cong C^n$ as $F$-algebras.

Now let $C'$ be another $H$-algebra, and let $\psi:C\to C'$ be an
$H$-algebra map.  It is immediate that the $G$-module map $\I_H^G(\psi)$ is
also an algebra homomorphism.  (Under the above identifications, it is just
$\psi\op\dots\op\psi:C^n\to (C')^n$.)  Thus, $\I_H^G$ is a functor.
\end{proof}

\begin{rems} \begin{mrk} \item If $H$ does not have finite index in $G$, then
$\I_H^G(B)$ is a nonunital $G$-algebra.  Indeed, the coinduced representation
is still a $G$-algebra, and $\I_H^G(B)$ is isomorphic to the nonunital
subalgebra of $FH$-maps which are finitely supported modulo $H$.
\item If $B$ is an interior $H$-algebra, i.e. $H$ acts on $B$ by inner
automorphisms, then there is another way of defining an induced $G$-algebra
originally introduced by Puig.  These two concepts are quite different.
Indeed, the underlying $G$-module in Puig's construction is not $\I_H^G(B)$,
but instead $\I_H^G(B)\ot_{FH}FG$.  The resulting $F$-algebra structure is
isomorphic to $M_n(B)$ instead of $B^n$ \cite[\S 16]{T}.
	     \end{mrk}
\end{rems}

It is easy to check that this functor satisfies the usual properties of
induction.
\begin{prop}\label{T:ind} Let $H$ be a subgroup of $G$ of finite index, and suppose that
$C$ and $C'$ are $H$-algebras.
\begin{enumerate} \item $\I_H^G(C\op C')\cong\I_H^G(C)\op \I_H^G(C')$ and
$\I_H^G(C\cap C')=\I_H^G(C)\cap \I_H^G(C')$ as $G$-algebras.
\item If $C$ is an $H$-subalgebra of $C'$, then $\I_H^G(C)$ is a
$G$-subalgebra of $\I_H^G(C')$, and $C=C'$ if and only if
$\I_H^G(C)=\I_H^G(C')$.
\item If $H\le K\le G$, then $\I_K^G(\I_H^K(C))=\I_H^G(C)$.
\end{enumerate}
\end{prop}

We now return to our construction of invariant subalgebras.  Suppose that
$V=\I_H^G(W)$, where $W$ is a $D$-module which is a projective
representation of $H$.  The cocycle defining $\rho_W$ is just the
restriction of $\a$ to $h\times H$.  It is automatic that $W$ is
irreducible.  Since $W$ is a direct summand of
$V_H\overset{\text{def}}{=}\text{Res}_H^G(V)$, $H$ must act on $W$ by
$D$-linear automorphisms; this means that $W$ is an $F^{\a}H\ot D$-module.
Note that the induced representation $V$ comes equipped with a
distinguished choice of $F^{\a}H$-submodule isomorphic to $W$ (namely
$\b{1}\ot W)$, and the invariant subalgebra we construct below depends on
this choice.  For ease of notation, we view $W$ as this fixed $H$-submodule
of $V$.  Let $T=\{g_1=1,g_2,\dots,g_l\}$ be a left transversal of $H$, and
set $W_i=\ov{g_i}\ot W$, a $D$-subspace of $V=F^{\a}G\ot_{F^{\a}H}W$.
Define a map $\Psi_{(H,W,T)}:\I_H^G(\E_D(W))\to A=\E_D(V)$ via the formula
$\Psi_{(H,W,T)}((g_i\ot f))(\ov{g_j}\ot w)=\d_{ij}\ov{g_i}\ot f(w)$ for
$f\in\E_D(W)$, $w\in W$ and extending by linearity.

\begin{lem} \label{T:psi} The map $\Psi_{(H,W)}=\Psi_{(H,W,T)}$ is independent of the choice of
transversal.  It is an injective $G$-algebra homomorphism whose
image is the block-diagonal subalgebra $\op_{i=1}^l \E_D(W_i)$.  In
particular, this subalgebra is $G$-invariant. 
\end{lem}
\begin{proof}  
Let $\Psi=\Psi_{(H,W,T)}$.  It is easy to see that $\Psi$ is an embedding of algebras with the
specified image, so we need only check that $\Psi$ is an intertwining map.
Fix $g\in G$.  There exists a permutation $\s=\s_g\in S_l$ and elements
$h_i\in H$ such that $gg_i=g_{\s(i)}h_{\s(i)}$ for all $i$. First note that
\begin{equation*}
\begin{split}
\Psi(g\cdot(g_i\ot f))(\ov{g_j}\ot w)=\Psi(gg_i\ot f)(\ov{g_j}\ot
w)&=\Psi(g_{\s(i)}\ot h_{\s(i)}\cdot f)(\ov{g_j}\ot w) \\
&=\d_{\s(i)j}\ov{g_j}\ot (h_j\cdot f)(w).
\end{split}
\end{equation*}  
On the other hand, a similar
calculation using the definition of multiplication in $F^{\a}G$ gives
\begin{equation*}
(g\cdot\Psi(g_i\ot f))(\ov{g_j}\ot w)=\d_{\s(i)j}\beta \ov{g_j}\ot
\ov{h_j}f(\ov{h_j}^{-1}w),
\end{equation*}
where
$$\beta=\a(g_j,h_j)^{-1}\a(h_j,h_j^{-1})\a(g,g^{-1}g_jh_j)\a(g^{-1}g_jh_j,h_j^{-1})^{-1}\a(g^{-1},g_j)\a(g,g^{-1})^{-1}.$$
Applying the cocycle condition and the fact that $\a(x,1)=1=\a(1,x)$ for
all $x\in G$, we get
\begin{equation*}
\begin{split}
\beta
&=\a(g_j,h_j)^{-1}\a(h_j,h_j^{-1})\a(g_jh_j,h_j^{-1})^{-1}\a(g,g^{-1}g_j)\a(g^{-1},g_j)\a(g,g^{-1})^{-1}\\
&=\a(g_j,h_j)^{-1}\a(h_j,h_j^{-1})\a(g_jh_j,h_j^{-1})^{-1}=1,
\end{split}
\end{equation*}
as desired.

The verification that $\Psi$ does not depend on the transversal is similar,
but easier.
\end{proof}

  Let $C$ be an invariant subalgebra of the $H$-algebra $\E_D(W)$.  It now
follows from Proposition \ref{T:ind} and the lemma that $\Psi(\I_H^G(C))$
is a $G$-invariant subalgebra of $A=\E_D(V)$.  More precisely,
\begin{prop} \label{T:fund} 
The map $C\mapsto \Theta_{(H,W,C)}\overset{\text{def}}{=}\Psi_{(H,W)}(\I_H^G(C))$ defines an injective lattice homomorphism from the
$H$-invariant subalgebras of $\E_D(W)$ to the $G$-invariant subalgebras of
$\op_{i=1}^l \E_D(W_i)\subset A$.
\end{prop}

It is not true that an invariant subalgebra of $A$ can be expressed
uniquely in terms of this construction if the initial data (namely $H$,
$W$, and $C$) are allowed to vary.  Indeed, conjugate data
(i.e. $gHg^{-1}$, $\b{g}W$, and $g\cdot C\subseteq\E_D(\b{g}W)$ for some
$g\in G$) produces the same invariant subalgebra.  However, we will see
below that uniqueness does hold if we restrict ourselves to conjugacy
classes of initial data with $C$ simple.

In order to show that this construction gives rise to all invariant
subalgebras, we first need to associate a transitive permutation representation
of $G$ to any invariant subalgebra $B$.  By \ref{T:ss}, we can write
$B=B_1\oplus\dots\oplus B_l$, where the $B_i$ are simple.  The restriction
of the $G$-action $\pi$ to $B$ gives rise to a permutation representation
of $G$ on the set of $B_i$'s because the algebra automorphism $\pi(g)$ must
permute the minimal two-sided ideals of $B$.  More explicitly, let
$X=\{e_1,\dots,e_l\}$ with $e_i=1_{B_i}$ be the set of central primitive
idempotents of $B$.  Since $e_i$ is the unique nonzero idempotent in the
center of $B_i$, it is clear that if $\pi(g)(B_i)=B_j$, then
$\pi(g)(e_i)=e_j$.  We thus obtain a homomorphism $\bar{\pi}_B:G\to S_l$,
where we have identified $S(X)$ with $S_l$ in the obvious way.
Note that $X$ is also the set of central primitive idempotents in $Z_A(B)$
and $Z(B)$.  Accordingly, dual pairs give rise to the same permutation
representation, as do any invariant subalgebras with the same center.

The permutation representation $\b{\pi}_B$ can also be defined in terms of
the $B$-isotypic components of $V$.  Recall that $V=\op_{i=1}^l V_i$ where
$V_i=V'_i$ is the isotypic $B$-submodule of $V$ corresponding to $B_i$.
Fix $g\in G$ and $v\in V_j$, and, write $\b{g}^{-1}(v)=\sum_{i=1}^l v_i^g$
with $v_i^g\in V_i$.  Note that $\b{g}^{-1}(g\cdot
e_i)(v)=e_i(\b{g}^{-1}(v)=v_i^g$.  But by definition, $g\cdot
e_i=e_{\bar{\pi}_B(g)(i)}$, giving $v_i^g=\b{g}^{-1}(g\cdot
e_i)(v)=\delta_{j,\bar{\pi}_B(g)(i)}\b{g}^{-1}(v)=\delta_{i,\bar{\pi}_B(g^{-1})(j)}\b{g}^{-1}(v)$.
This implies that $\b{g}^{-1}(V_j)\subseteq V_{\bar{\pi}_B(g^{-1})(j)}$ for
all $j$.  Applying this to $g^{-1}$ (or using the fact that $\r(g)^{-1}$ is
surjective) gives the reverse inclusion.  Thus, $G$ permutes the $V_i$'s,
and this permutation is just $\bar{\pi}_B$.

\begin{prop}  The permutation representation $\bar{\pi}_B$  is  transitive.
\end{prop}

\begin{proof}  

Let $U$ be a simple $B$-submodule of $V$ isomorphic to a minimal left ideal
of $B_1$.  By definition, $e_1$ is the identity map on $U$.  Let $e_i$ be
any central primitive idempotent, and choose $g\in G$ such that $\b{g}U$ is a
simple $B_i$ module.  For all $u\in U$, we have $(g\cdot
e_1)(\b{g}u)=\b{g}(e_1(\b{g}^{-1}(\b{g}u)=\b{g}(e_1(u))=\b{g}(u)$.  Since
$e_i$ is the unique central primitive idempotent acting as the identity on
$\b{g}U$, this implies that $g\cdot e_1=e_i$.
\end{proof}

If $G$ acts on $B$ by inner automorphisms, then the $G$-action preserves
the simple components of $B$.  We thus obtain the useful corollary:
\begin{cor}  If $G$ acts on the invariant subalgebra $B$ by inner
automorphisms, then $B$ is simple.
\end{cor}

Let $H_i=\{g\in G\mid g\cdot e_i=e_i\}$ be the inertia subgroup of $e_i$.
Note that it has finite index $l$ in $G$.  It is immediate that $V_i$ is an
$F^{\a}H_i\ot D$ submodule of $V$, and the
transitivity of $\b{\pi}_B$ implies that $V=\I_{H_i}^G(V_i)$, i.e. $V$ is
isomorphic to the induced representation and has distinguished
$H_i$-submodule $V_i$.  Moreover, $V_i$ is an ($F$-)irreducible projective
representation of $H_i$ because if $M$ were a proper subrepresentation,
then $\I_{H_i}^G(M)$ would be a proper $G$-submodule of $V$, contradicting
the irreducibility of $V$.  The algebra $B_i$ is a simple $H_i$-subalgebra of
$\E_D(V_i)$, and we are precisely in the situation of the fundamental
construction.  The uniqueness part of Proposition \ref{T:algind}
shows that $B=\Theta_{(H_i,V_i,B_i)}$.  We have thus realized $B$ in $l$ different ways, all of
which have conjugate initial data.

Now suppose that $B=\Theta_{(H,W,C)}$.  By definition, $W$ is the isotypic
$B$-submodule corresponding to the simple component $C$ (i.e. $1\ot C$) of $B=\I_H^G(C)$,
implying that $W=V_j$ and $C=B_j$ for some $j$.  Also, $H$ is the
stabilizer of $B_j$, so in fact $H=H_j$. 

Let $\eusm D$ be the set of equivalence classes of triples $(H,W,C)$ where $H$ is a subgroup of finite
index in $G$, $W$ is an $F^{\a}H\ot D$ submodule of $V$  such that $V\cong\I_H^G(W)$, and $C$ is an invariant
subalgebra of the $H$-algebra $\E_D(W)$.    Also, let $\eusm
D_{(H,W)}\subset\eusm D$ be the subset of classes with a representative of
the form $(H,W,C)$.

\begin{thm} \label{T:main} Let $A=\E_D(V)$ be a $G$-simple central simple algebra.  The
map $(H,W,C)\mapsto\Theta_{(H,W,C)}$
gives a  bijective
correspondence between $\eusm D$ and the set of unital $G$-invariant
subalgebras of $A$.    This bijection preserves dual pairs and centers; if
$B=\Theta_{(H,W,C)}$, then $Z_A(B)=\Theta_{(H,W,Z_{\E_D(W)}(C))}$ and
$Z(B)=\Theta_{(H,W,Z(C))}$.  Similarly, $Z_0(B)$ (the $F$-linear span of
the Wedderburn components of $B$) is just $\Theta_{(H,W,Z_0(C))}=\Theta_{(H,W,F1_{\E_D(W)})}$.
Furthermore, the image of  $\eusm
D_{(H,W)}$ under the correspondence is precisely the set of  invariant
subalgebras $B$ with $Z_0(B)=\Theta_{(H,W,F1_{\E_D(W)})}$. 
\end{thm}
\begin{proof}  
We have already shown that there is a bijection between invariant
subalgebras and triples $(H,W,C)$ where $W$ is an  $F^{\a}H$-submodule of
$V$ such that the obvious map $W\to\b{1}\ot W$ extends to an isomorphism
$V\cong\I_H^G(W)$.  This amounts to saying that $V$ is the internal direct
sum of the translates $\ov{g_i}W$ (and so $V$ can be viewed as  equal and not
just isomorphic to $\I_H^G(W)$).  The following lemma shows  that any subrepresentation of $V_H$
isomorphic to $W$ satisfies this condition.  

\begin{lem}
Suppose that $V=\I_H^G(W)$ with $V$ irreducible.  Then if $W'$ is any
subrepresentation of $V_H$ isomorphic to $W$, $V$ is the internal direct
sum of the $\ov{g_i}W'$'s.
\end{lem}
\begin{proof}
By Frobenius reciprocity, there is a linear isomorphism
$\H_{F^{\a}H}(W,V_H)\cong\H_{F^{\a}G}(V,V)$ given by $f\mapsto
\hat{f}$, with $\hat{f}(\ov{g_i}\ot w)=\ov{g_i}f(w)$.  Let
$f:W\to V_H$ be an $H$-map with image $W'$.     Since $V$ is
irreducible, $\hat{f}$ is an isomorphism.  Accordingly, $V$ is the direct
sum of the distinct $G$-translates of $f(W)=W'$.
\end{proof}

It only remains to prove the last three statements.  We have shown
that as an $F$-algebra, $\Theta_{(H,W,C)}$ is just $C^{[G:H]}$ embedded in
the block diagonal subalgebra $\op_{i=1}^l \E_D(W_i)\subset A$.  Since
taking finite direct sums commutes with taking dual
pairs, centers, and $Z_0$, the result follows.
\end{proof}

\begin{rems}\begin{mrk} 
\item Since an invariant subalgebra $B$ can always be expressed trivially
as $\Theta_{(G,V,B)}$, it is clear that a nonsimple $B$ can arise from
nonconjugate initial data.  The class  in  $\eusm D$
corresponding to $B$ consists of the triples with minimal $H$ (or $W$ or $C$). 
\item Let $F$ be an infinite field.  If $V\cong\I_H^G(W)$ and $V_H$ does
not have a unique subrepresentation isomorphic to $W$, then $A$ has an
infinite number of invariant subalgebras.  Indeed, in this case, the
$W$-isotypic submodule of $V_H$ is a direct sum of $t\ge 2$ submodules
isomorphic to $W$, so there are an infinite number of submodules $W'$
isomorphic to $W$.  At most $[G:H]$ of these submodules can be conjugate,
and each class gives rise to a distinct invariant subalgebra
$\Theta_{(H,W',F)}$.
	    \end{mrk}

\end{rems}

Before proceeding, we give two examples in the case $A=\E_F(V)$.
\begin{egs}\begin{mrk}\item
Let $V$ be primitive, i.e. suppose that $V$ is not induced from any proper
subgroup.  Then all invariant subalgebras of $A$ are simple.
\item The theorem shows that $V$ is a monomial representation, i.e. it is
induced from a linear character, if and only if $\E_F(V)$ has a
$G$-invariant split Cartan subalgebra $\frak h$.  Indeed, this can be shown
directly.  By choosing an appropriate basis for $V$, we can view $\frak h$
as the subalgebra of diagonal matrices in $M_n(F)$.  Note that for $\frak
h$ to be $G$-invariant means precisely that its normalizer $N(\frak h)$
contains $\r(G)$.  But $N(\frak h)$ is the set of monomial matrices, and it
is well known that $V$ is monomial if and only if $\r(G)$ consists of
monomial matrices with respect to some basis for $V$. \cite[p.67]{I}.
\end{mrk}
\end{egs}

The correspondence in this theorem becomes much simpler when $V$ has nice
rationality properties.  Recall that a projective $F$-representation $V$ is
called absolutely irreducible if $V_E=V\otimes E$ is an irreducible
projective $E$-representation for every algebraic extension $E$ of
$F$. Equivalently, the division algebra
$\E_G(V)\overset{\text{def}}{=}\E_{F^{\a}G}(V)$ is just the ground field
$F$.  Note that if $F$ is algebraically closed, then all irreducible
representations are absolutely irreducible.

\begin{lem}  Let $A$ be $G$-simple.  If $K=\E_G(V)$, then
$D=\E_A(V)\subseteq K$.  In particular, if $V$ is absolutely irreducible, then
$D=F$ and $A=\E_F(V)$.
\end{lem}
\begin{proof}  
Choose $d\in D$.  Then we have $d(\r(g)v)=\r(g)(dv)$ for $g\in G$, $v\in
V$, since $\r(g)\in A$.  Hence, $d\in K$.
\end{proof}
If $V$ is absolutely irreducible, we call such $A=\E_F(V)$ absolutely $G$-simple.

Now suppose that $H$ is a subgroup of finite index and $W$ is an
 (irreducible) $F^{\a}H$-module such that $V\cong\I_H^G(W)$.  Here, we are
 not viewing $W$ as a specific subspace of $V$.  If $V$ is absolutely
 irreducible, then $\H_{F^{\a}G}(\I_H^G(W),V)$ is one-dimensional.  By
 Frobenius reciprocity, the same is true for $\H_{F^{\a}H}(W,V_H)$. This
 implies that there is a unique subrepresentation of $V_H$ isomorphic to
 $W$, since otherwise there would be linearly independent $H$-maps $W\to
 V_H$.  Similarly, we must have $\E_H(W)=F$.  Summing up:

\begin{prop} \label{T:abs} Let $V$ be absolutely irreducible, and suppose that
$V\cong\I_H^G(W)$ where $H$ is a subgroup of finite index and $W$ is an
irreducible $F^{\a}H$-module.  Then there is a unique subrepresentation of
$V_H$ isomorphic to $W$.  Moreover, $W$ is absolutely irreducible.
\end{prop}

Let $\tilde{\eusm D}$ be the set of conjugacy classes of triples where $W$
is only defined up to isomorphism, i.e. $W$ is no longer viewed as a
specific subspace of $V$.  In other words, $\tilde{\eusm D}$ consists of
the classes of $\eusm D$ modulo $H$-isomorphism of the second variable.
It is clear that triples in $\eusm D$ representing the same class in $\tD$
give rise to invariant subalgebras that are isomorphic as $G$-algebras.  If
$V$ is absolutely irreducible, the previous proposition shows that the
projection $\D\to\tD$ is a bijection.  Accordingly, we get the first
statement of the corollary:
\begin{cor} \label{T:tD}
Let $A=\E_F(V)$ be absolutely $G$-simple.  The map $(H,W,C)\mapsto\Theta_{(H,W,C)}$ gives a bijective
correspondence between $\tD$ and the set of unital $G$-invariant
subalgebras of $A$.  In addition, $\Theta_{(H,W,C)}$ is separable;
equivalently, $Z(C)$ is a separable field extension of $F$.
\end{cor}
\begin{proof}
Write $B=\Theta_{(H,W,C)}$.  Extending scalars to the algebraic extension
$E$ gives the invariant subalgebra $B_E$ of the central simple $E$-algebra
$A_E\cong\E_{D_E}(V_E)$.  Since $V_E$ is irreducible, Proposition
\ref{T:ss} applies, showing that $B_E$ is semisimple.  Thus, $B$ is separable.
\end{proof}

We are now ready to make the correspondence in Theorem \ref{T:main}
entirely explicit when $F$ is algebraically closed.  We start by
classifying invariant central simple subalgebras of any $G$-simple $A$.

Let $B$ be a simple subalgebra of $A=\E_D(V)$ with simple $B$-module $W'$
and simple $B\ot D$-module $W=W'\ot D$.  The $B\ot D$-module $V$ is
isotypic, say $V\cong mW$.  Let $L=\E_{B}(W')=\E_{B\ot D}(W)$ and set
$U=(L^{op})^m$.  We obtain the factorization $V\cong W\ot_{L^{op}} U\cong
(W'\ot_{L^{op}} U)\ot_F D$.  As shown in the proof of Theorem \ref{T:DCT},
$Z_A(B)=\E_{L^{op}}(U)$; also, $B\cong \E_{L}(W')\cong\E_{L\ot
D}(W)$.  In addition, any dual pair of simple subalgebras arises in this way.

\begin{prop} Let $A=\E_D(V)$ be a central simple algebra.
If $V\cong W\ot_{L^{op}} U\cong
(W'\ot_{L^{op}} U)\ot_F D$ with $W'$ an $L$-module, $U$
an $L^{op}$-module, and  $W=W'\ot D$ an $L\ot D$-module, then
$\E_{L}(W')$ and $\E_{L^{op}}(U)$ is a dual pair of simple
subalgebras.  Conversely, any dual pair of simple subalgebras comes from
such a factorization.  In addition, the subalgebras are central simple if
and only if $L$ is a central division algebra.
\end{prop}

Using this result, we can classify invariant central simple subalgebras.
Let $L$ be a central division algebra, and let $W'$ and $U$ be $L$ and
$L^{op}$ modules respectively which are projective representations given by
$G\overset{\r_{W'}}{\to}\E_L(W')^{\times}$ and
$G\overset{\r_{U}}{\to}\E_{L^{op}}(U)^{\times}$.  Set $V=(W'\ot_{L^{op}}
U)\ot_F D$, and let $\tau$ denote the canonical isomorphism $\E_L(W')\ot_F
\E_{L^{op}}(U)\overset{\tau}{\to}\E_D(V)$ given by $\tau(f_1\ot f_2)(w'\ot
u\ot d)=f_1(w)\ot f_2(u)\ot d$.  Then $\r_V:G\to\E_D(V)^{\times}$ defined
by $\r_V(g)=\tau(\r_{W'}(g)\ot\r_U(g))$ makes $V$ into a projective
representation.  It is easy to check that $\tau$ becomes a $G$-algebra
isomorphism. If $\r_{W'}$ and $\r_U$ are twisted by (one-dimensional)
projective characters, then the new $G$-action on $V$ is projectively
equivalent to the old one.

Conversely, suppose that $V$ is a projective representation, and $\E_L(W')$
and $\E_{L^{op}}(U)$ are invariant.  The map $\tau$ is thus a $G$-algebra
isomorphism.  By Proposition \ref{T:SN}, the $G$-actions
on these subalgebras come from projective representations $(W',\r_{W'})$ and
$(U, \r_U)$.  Hence, $\tau^{-1}(\r_{W'}\ot \r_U)$ and $\r_V$ define the
same $G$-algebra structure on $\E_D(V)$, implying that they are
projectively equivalent, i.e. differ by a projective character.  Modifying
$\r_U$ by this twist, we get $\r_V=\tau(\r_{W'}\ot\r_U)$.
It is obvious that if $V$ is irreducible, then both $W'$ and $U$ must be as
well.  This proves the
following theorem:

\begin{thm} 
Let $A=\E_D(V)$ be $G$-simple.  Suppose that $V\cong(W'\ot_{L^{op}} U)\ot_F
D$ is a factorization such that $L$ is a central division algebra and $W'$
and $U$ are (irreducible) projective representations of $G$ (via $L$ and
$L^{op}$ linear automorphisms respectively). Then
$A\cong\E_{L}(W')\ot\E_{L^{op}}(U)$ as $G$-algebras and the images of the
two factors in $A$ are a dual pair of invariant central simple subalgebras.
Conversely, any such dual pair arises in this way. 
\end{thm}
\begin{rem} 
If $D=F$, invariant central simple subalgebras come from
expressing $V$ as the tensor product of projective representations.  
In general, finding all (or even some) factorizations for a given $V$ is a
difficult problem.  See for example \cite{St}.
\end{rem}

We can say more when  $V$ is absolutely irreducible.  Recall that in this
case, $D=F$ and $W=W'\ot_F D=W'$.  Since $\E_{F^{\a}G}(V)=F$, any two $G$-maps $W\ot_{L^{op}}
U\overset{\sim}{\to}V$ are scalar multiples of each other and thus give the
same dual pair of invariant central simple subalgebras.  Thus, the specific
factorization does not matter.

\begin{cor} \label{T:cs}
If $A$ is absolutely $G$-simple, then there is a one-to-one correspondence
between pairs of irreducible projective representations $(W,U)$  modulo projective
equivalence such that
$V\cong(W\ot_{L^{op}} U)$ and dual pairs of invariant central simple
subalgebras.
\end{cor}

We now describe the index set for the classification of invariant
subalgebras in the algebraically closed case.  Let $\eusm E'$ be the set of
quadruples $(H,W,W_1,W_2)$ where $H$ is a subgroup of $G$ of finite index,
$W$ is an irreducible projective representation of $H$ such that
$V\cong\I_H^G(W)$, and $W_1$, and $W_2$ are irreducible projective
representations of $H$ such that $W\cong W_1\ot_F W_2$.  We then let $\eusm
E$ be the set of equivalence classes of $\eusm E'$ where two quadruples
$(H,W,W_1,W_2)$ and $(H',W',W'_1,W'_2)$ are equivalent if there exists
$g\in G$ such that $H'=H^g$, $W'=W^g$, and $W'_i$ is projectively
equivalent to $W_i^g$.  We
let $\eusm E_{(H,W)}\subset \eusm E$ be the subset of classes with a
representative of the form $(H,W,W_1,W_2)$.  In addition, we denote by
$C(W_1,W_2)$ the image of $\E_F(W_1)\ot 1$ under the isomorphism
$\E_F(W_1)\ot \E_F(W_2)\to\E_F(W)$.  The trivial factorizations give
$C(F,W)=F$ and $C(W,F)=\E_F(W)$.

\begin{thm}  \label{T:eE} Let $F$ be algebraically closed and $A=\E_F(V)$ a $G$-simple
algebra.  Then the map $(H,W,W_1,W_2)\mapsto\Theta_{(H,W,C(W_1,W_2))}$
gives a bijective correspondence between $\eusm E$ and the set of invariant
subalgebras of $A$.  Moreover, the duality on invariant subalgebras is
given by interchanging the $W_i$'s,
i.e. $Z_A(\Theta_{(H,W,C(W_1,W_2))})=\Theta_{(H,W,C(W_2,W_1))}$.  The image
of $\eusm E_{(H,W)}$ under the correspondence is precisely the set of
invariant subalgebras $B$ with center $\Theta_{(H,W,C(F,W))}$.
\end{thm}
\begin{proof}   Recall that $\tD$ is the set
of classes of triples $(H,W,C)$ where $H$ and $W$ are defined as in $\eusm
E$ and $C$ is a (central) simple subalgebra of $\E_F(W)$ (using the fact
that $F$ is algebraically closed).  Since $V$ is absolutely irreducible,
Corollary \ref{T:tD} shows that invariant subalgebras are parameterized by
this set.  Applying Corollary \ref{T:cs}, we see that the map
$(H,W,W_1,W_2)\mapsto (H,W,C(W_1,W_2))$ induces a bijection $\eusm
E\to\tD$, and we obtain the desired correspondence.  Since
$Z_{\E_F(W)}(C(W_1,W_2))=C(W_2,W_1)$ and $Z(C(W_1,W_2))=C(F,W)$ , the last statements follow from Theorem
\ref{T:main}.
\end{proof}
\begin{rem}  Note that the cocycle $\a$ does not determine the cocycles
defined by $\r_{W_1}$ and $\r_{W_2}$.  In particular, even if $V$ is a
linear representation, it is not possible to avoid considering projective
representations when studying invariant subalgebras of $\E_F(V)$.
\end{rem}

It is convenient to reformulate this correspondence in terms of covering
groups.  Recall that $\td{G}$ is an $F^*$-generalized covering (or
representation) group for $G$ if it is a central extension of $G$
satisfying the projective lifting property for projective representations
over $F$.  It is known that $F^*$-generalized covering groups always
exist.  If $F$ is algebraically closed and $G$ is finite, then we can
choose $\td G$ finite of order $|G||H^2(G,F^*)|$; such a group is called an
$F^*$-covering group for $G$ \cite{BT}. 

We now assume that $F$ is algebraically closed (so $D$ and $L$ are just $F$
and $W=W'$).  Suppose that the projective representation $V$ factors as
$V\cong W\ot_F U$.  Choose a linear representation $(V,\td{\r}_V)$ of
$\td{G}$ lifting $\r_V$ and similarly for $W$ and $U$.  A priori, $V$ is
only projectively equivalent to $W\ot U$ over $\td{G}$.  However, if $V_1$
and $V_2$ are linear representations which are projectively equivalent,
then $V_1\cong V_2\ot \lambda$, where $\lambda$ is a linear character.
Thus, by choosing a different lift for $\r_W$, we obtain linear
representations of $\td{G}$ such that $V\cong W\ot_F U$ as
$\td{G}$-modules.  On the other hand, it is obvious that any such
factorization gives an isomorphism of projective representations for $G$.

This allows us to redefine $\eusm E_{(H,W)}$.  Let $\td H$ be a generalized
covering group for $H$, and fix a lift of $W$ to a linear representation of
$\td H$.  If $(W_1,W_2)$ and
$(W'_1,W'_2)$ are two pairs of linear representations of $\td H$ satisfying
$W\cong W_1\ot W_2\cong W'_1\ot W'_2$, we say they are equivalent
if for some linear character $\lambda$ of $\td H$, $W'_1\cong W_1\ot\lambda$
and $W'_1\cong W_1\ot\lambda^{-1}$.  Denote the set of such classes by
$\eusm F_{(H,W)}$.  The previous observations give the
following result.
 
\begin{lem}  There is a natural bijection between $\eusm E_{(H,W)}$ and
$\eusm F_{(H,W)}$. 
\end{lem}

Let $Y$ be a complete set of representatives of the conjugacy classes of
pairs $(H,W)$.  Then the $\eusm E_y$'s partition $\eusm E$.  Set $\eusm
F=\coprod_{y\in Y}\eusm F_y$.  We can now rewrite Theorem \ref{T:eE}.
\begin{thm} \label{T:eF}
 Let $F$ be algebraically closed and $A=\E_F(V)$ a $G$-simple
algebra.   Then the map $(H,W,W_1,W_2)\mapsto\Theta_{(H,W,C(W_1,W_2))}$
gives a bijective correspondence between $\eusm F$ and the set of invariant
subalgebras of $A$.  Duals and centers of invariant subalgebras are given
by the same formulas as before.
\end{thm}

It is possible to avoid all explicit mention of projective representations
in classifying invariant subalgebras.  In order to do this, choose a
generalized covering group $\td G$ of $G$ and fix a lift of $V$ to a
representation of $\td G$.  Since the $G$ and $\td G$ invariant subspaces
of $A$ are the same, we can apply the above procedure to the $\td G$-simple
algebra $A$.  Note that this will require choosing a generalized covering
group ${\td{\td G}}$ of $\td G$!

If $F$ is not algebraically closed, it is not true in general that a simple
$G$-algebra $A$ will have a finite number of invariant subalgebras, even when
$G$ is finite.  We have already seen a way that finiteness can fail if $F$
is infinite and $V$ is not absolutely irreducible.  Namely, if $V\cong
\I_H^G(W)$ and $V_H$ does not have a unique subrepresentation isomorphic to
$W$, then for any simple $H$-invariant $C\subset \E_D(W)$, the set
$\{\Theta_{(H,W',C)}\mid W'\subset V$, $W'\cong W\}$ will be infinite.
Note that these subalgebras are all nonsimple.

Furthermore, the set of invariant subalgebras can be infinite even when $V$ is primitive.
  Indeed, we have
the proposition:
\begin{prop}  Let $A=\E_F(V)$ where $V$ is an irreducible projective 
representation of $G$, and suppose that the division algebra
$\E_G(V)$ is not a field.  Then $\eusm D_{(G,V)}$ is infinite, i.e $\E_F(V)$ has an infinite number of
simple invariant subalgebras.
\end{prop}
\begin{proof} Note that any subalgebra of $\E_G(V)=(\E_F(V))^G$ is
$G$-invariant, so the following lemma gives the result.
\end{proof}

\begin{lem} Let $D$ be a noncommutative central $F$-division
algebra.  Then $D$ contains an infinite number of distinct
subfields.
\end{lem}
\begin{proof}  Choose noncommuting elements $u,v\in D$, and consider
the subfields $F_a=F(u+av)$ for $a\in F$.  Wedderburn's theorem on
finite division rings shows that the field $F$ is infinite, so it
suffices to show that $F_a=F_b$ if and only if $a=b$.  If
$F_a=F_b$, then $u+av$ and $u+bv$ commute, implying that
$auv+bvu=buv+avu$.  If $uv$ and $vu$ are linearly independent over
$F$, it is immediate that $a=b$.  Otherwise, $vu=cuv$ for some
$c\in F$, giving $(a+bc)uv=(b+ac)uv$ and $(c-1)(a-b)=0$.  Since
$c\ne 1$, $a=b$.
\end{proof}
 
However, these pathologies cannot occur when $F$ is algebraically closed.
\begin{thm}
Let $F$ be algebraically closed, $G$ a finite group, and $A=\E_F(V)$ a $G$-simple
algebra.  Then $A$ has a finite number of invariant subalgebras.
\end{thm}
\begin{proof}
 Replacing $G$ by a covering group (which is also
finite), we can assume without loss of generality that $V$ is a linear
representation of $G$.  Since the set of invariant subalgebras and $\eusm F=\coprod_{y\in Y}\eusm F_y$
have the same cardinality (using the notation of Theorem \ref{T:eF}), it
suffices to show that $Y$ and the $\eusm F_y$'s are finite.  A theorem of
Berman and Witt shows that for arbitrary $F$, the number of nonisomorphic
irreducible $F$-representations of a finite group is finite \cite[Theorem
17.5.3]{K}.  The set $Y$ is finite because it is contained in the set of all pairs
$(H,W)$ where $H$ is a subgroup of $G$ and $W$ is an isomorphism class of
irreducible $FH$-modules.   Also, $\eusm F_{(H_y,W_y)}$ is finite, since it
is smaller than the set of arbitrary pairs of isomorphism classes of
irreducible $F\td H$-modules, where $\td H$ is a covering group for $H$.
\end{proof}

We conclude this section with an application to nonunital invariant
subalgebras.
\begin{prop} \label{T:nu}
Let $F$ be an algebraically closed field and $V$ an irreducible primitive
projective representation of $G$.  Then $\{0\}$ is the only nonunital invariant
subalgebra of  $A=\E_F(V)$.  Equivalently, any nonzero subrepresentation of $A$
closed under multiplication must contain the identity.
\end{prop}
\begin{proof}
We begin with a lemma.
\begin{lem} Let $F$ be an algebraically closed field.  For $t\ge 2$, the matrix algebra $M_t(F)$ has no nonunital subalgebras of codimension one.
\end{lem}
\begin{proof}
Suppose that $Q$ is a nonunital subalgebra of codimension one.  First note
that any element of $Q$ must be singular.  To see this, take $q\in Q$
invertible, so that $\det q\ne 0$.  It is a well-known corollary of the
Cayley-Hamilton theorem that $q^{-1}$ can be expressed as a polynomial in
$q$, so $q^{-1}\in Q$.  This implies that $Q$ contains the identity, a
contradiction.  Thus, $Q\subseteq V(\det)$, the hypersurface of $M_t(F)$ cut
out by the determinant.  But $Q$ is also a codimension one linear
subvariety, so $Q=V(f)$ for some homogeneous degree one polynomial $f$.  As
a result, $f$ divides $\det$, and this cannot be true, since the determinant is an
irreducible polynomial of  degree $t$.
\end{proof}

Now, let $Q$ be an nonunital invariant subalgebra.  Then $Q'=Q+F1_A$ is a
unital invariant subalgebra.  We know from the first example after Theorem
\ref{T:main} that $Q'$ is simple, hence isomorphic to $M_t(F)$ for some
$t\ge 1$.  If $t=1$, then $Q=\{0\}$.  Applying the lemma finishes the
proof.
\end{proof}

\section{Invariant subalgebras for topological and Lie groups}

In this section, we classify invariant subalgebras in the case where
$V$ is a continuous irreducible complex projective representation of a
compact connected Lie group.  For the moment, we consider a more general
situation.  Suppose that $G$ is a topological group,
$A=\E_D(V)$ is a $G$-simple algebra endowed with a $T_1$ topology, and $G$ acts continuously on
$A$.  For example, the topology on $A$ could come from $F$ having the
structure of a $T_1$ topological field or $\E_F(V)$ could be given the
Zariski topology.  So far, this setting includes every abstract group $G$
and $G$-algebra considered in the previous section by giving $G$ and $A$
 the discrete topology.  In order to avoid this type of triviality, we
further assume that the connected component of the identity $G^o$ (a closed
normal subgroup) acts irreducibly on $V$.  We call such an algebra
topologically $G$-simple.

\begin{prop} Every invariant subalgebra of a topologically $G$-simple
algebra $A$ is simple.
\end{prop}
\begin{proof} 
A $G$-invariant algebra is also $G^o$-invariant, so it suffices to assume
that $G$ is connected.  Let $X$ be the set of central primitive
idempotents of an invariant subalgebra $B$.  The transitivity of $\pi_b$
shows that $X$ is connected.   However, since $A$ is
$T_1$, $X$ is  discrete.  This
implies that $X$ is a singleton, i.e. $B$ is simple.
\end{proof}

If we further assume that $F$ is algebraically closed, Theorem \ref{T:cs}
now applies to give a classification of the invariant subalgebras of
$A=\E_F(V)$ in terms of factorizations $V\cong W_1\ot W_2$ modulo
projective equivalence.

We now assume that $F=\C$ and $G$ is a compact Lie group.  Note that a
continuous homomorphism $G\to \A_{F-\text{alg}}(A)\subset GL(A)$ is a
continuous homomorphism $G\to PGL(V)$.  Thus, if $A$ is a continuous
$G$-algebra, then $V$ is a continuous projective representation.

\begin{lem}  Suppose that $G$ is a simple compact connected Lie group, and let
$V(\l)$ and $V(\mu)$ be irreducible representations with highest weights
$\l$ and $\mu$.  Then $V(\l)\ot V(\mu)$ is irreducible if and only if $\l$
or $\mu$ is $0$.
\end{lem}
\begin{proof}
Since $V(\l+\mu)$ is a component of $V(\l)\ot V(\mu)$, it suffices to
compare the dimension of these representations.  The Weyl dimension formula
states that $$\dim V(\l)=\prod_{\a\in
R^+}\frac{\langle\a,\l+\r\rangle}{\langle\a,\r\rangle},$$ where $R^+$ is
the set of positive roots, $\r$ is half the sum of the positive roots, and
$\langle\ ,\ \rangle$ is the Killing form.  The equation
$\langle\a,\l+\mu+\r\rangle\langle\a,\r\rangle+\langle\a,\l\rangle\langle\a,\mu\rangle=\langle\a,\l+\r\rangle\langle\a,\mu+\r\rangle$
shows that
$$\frac{\langle\a,\l+\mu+\r\rangle}{\langle\a,\r\rangle}\le\frac{\langle\a,\l+\r\rangle}{\langle\a,\r\rangle}\frac{\langle\a,\mu+\r\rangle}{\langle\a,\r\rangle},$$
with equality if and only if $\langle\a,\l\rangle\langle\a,\mu\rangle=0$.
If $\beta$ is the highest root, then $\langle\beta, \nu\rangle>0$ for any
nonzero dominant weight $\nu$.  Multiplying over all positive roots, it
follows easily that $\dim V(\l+\mu)<\dim V(\l)\dim V(\mu)$ if and only if
both $\l$ and $\mu$ are nonzero.
\end{proof}

Let $G$ be a compact connected Lie group.  It is well known that the
universal covering group of $G$ is of the form $\tilde G=
G_1\times\dots\times G_s\times\R^n$, where each $G_i$ is a simple, simply
connected, compact Lie group.  Let $V$ be an irreducible projective
representation of $G$.  Then $V$ can be lifted to an irreducible
representation of $\tilde G$, which can be expressed as $V_1\ot\dots\ot
V_s\ot L$, where $V_i$ is a complex irreducible representation of $G_i$ and
$L$ is a character of $\R^n$.  This means that $V$ is projectively
equivalent to $\tilde V=V_1\ot\dots\ot V_s$.  Moreover, simple Lie groups
have no nontrivial characters, so projective and linear equivalence are the
same for representations of $G_1\times\dots\times G_s$.  The lemma shows
that any factorization of $\tilde V=W\ot W'$ into the tensor product of two
representations of $\tilde G$ must have $W$ and $W'$ as complementary
partial products of $V_1\ot\dots\ot V_s$.  More precisely, let $I=\{i\mid
V_i\ne\C\}$ and take $J\subset I$.  Set $W_J=\bigotimes_{i=1}^s W_{Ji}$ and
$W'_J=\bigotimes_{i=1}^s W'_{Ji}$, where $W_{Ji}$ is $V_i$ if $i\in J$ and
$\C$ otherwise and $W'_{Ji}$ is $V_i$ if $i\ne J$ and $\C$ otherwise.  We
get a factorization $\tilde V=W_J\ot W'_J$, and $J\mapsto W_J$ gives a
one-to-one correspondence between the subsets of $I$ and the factors of
$\tilde V$.  This observation combined with Theorem \ref{T:cs} proves the
following theorem due to Etingof:
\begin{thm}
Let $G$ be a compact connected Lie group, and let $A=\E_{\C}(V)$ where $V$
is an irreducible projective representation of $G$ projectively equivalent to
$V_1\ot\dots\ot V_s$.  Then there is a bijective correspondence between
$\eusm P(I)$, the power set of $I=\{i\mid V_i\ne\C\}$, and the set of invariant
subalgebras of $A$, given by $J\mapsto \E_{\C}(W_J)$.  Moreover, the
duality operator corresponds to taking complements in $\eusm P(I)$, i.e. it
is given by $\E_{\C}(W_J)\mapsto\E_{\C}(W_{I-J})$.
\end{thm}

By Theorem \ref{T:cs}, there are no nontrivial invariant subalgebras not
containing $1_A$, so we obtain the corollary:
\begin{cor}
There are exactly $2^{|I|}+1$ subrepresentations of $\E_{\C}(V)$ which are
closed under matrix multiplication: $2^{|I|}$ unital subalgebras and
$\{0\}$.
\end{cor}

In particular, if $G$ is a simple compact connected Lie group, then no
topologically $G$-simple algebra has any nontrivial invariant subalgebras.
It would be interesting to find classes of finite group satisfying this
property and to find a group-theoretic characterization of such groups.  It
is not true that finite simple groups have this property.  In the notation of the
atlas of finite groups, $U_4(2)$ has irreducible representations $\chi_3$
and $\chi_4$ of dimensions five and six respectively such that
$\chi_3\ot\chi_4\cong \chi_{12}$ is also irreducible \cite{C}.

\end{document}